\documentstyle{amsppt}
\magnification=\magstep1
\NoRunningHeads
\NoBlackBoxes


%
%
\pageheight{9truein}
\pagewidth{6.5truein}
\TagsOnRight
%
\define\lra#1{{\left\langle#1\right\rangle}}
\define\lrbrace#1{{\left\{#1\right\}}}
\define\xp#1{{\phantom{#1}}}
\define\({\left(}
\define\){\right)}
\define\[{\left[}
\define\]{\right]}
\define\hra{\,{\hookrightarrow}\,}


\let\th=\thinspace

\let\dsp=\displaystyle
\let\germa=I
%
%
%
%
%
%
%
\define\ocal{{\Cal O}}
\define\bx{{\bold x}}
%

\define\Id{\operatorname{Id}}


\define\coker{\operatorname{coker}}
\define\Coker{\operatorname{Coker}}

\define\grade{\operatorname{grade}}

\define\Syz{\operatorname{Syz}}
\define\height{\operatorname{height}}

\define\Ext{\operatorname{Ext}}
\define\Hom{\operatorname{Hom}}
\define\im{\operatorname{Im}}

\define\Tor{\operatorname{Tor}}

\topmatter
\title The Monomial Conjecture and Order Ideals\endtitle
\author S.~P.\ Dutta\\
\endauthor
\address\newline
Department of Mathematics\newline
University of Illinois\newline
1409 West Green Street\newline
Urbana, IL 61801\newline
U.S.A.\newline
{\eightpoint e-mail: dutta\@math.uiuc.edu}
\endaddress
%
\thanks
\hbox{\hskip-12pt}AMS Subject Classification: Primary 13D02, 13D22;
secondary 13C15, 13D25, 13H05\hfill\hfill\newline
Key words and phrases: Regular local ring, order ideal, grade, syzygy,
canonical module.\hfill\hfill 
\endthanks

\abstract In this article first we prove that a special case of the order ideal conjecture,
originating from the work of Evans and Griffith in equicharacteristic, implies the monomial
conjecture due to M. Hochster. We derive a necessary and sufficient condition for the validity of this
special case in terms certain syzygis of canonical modules of normal domains possessing free summands.
We also prove some special cases of this observation.

\endabstract

\endtopmatter

\document


\bigskip

\baselineskip17pt


The main focus of this paper is to establish a relation between the
monomial conjecture due to Hochster \cite{Ho1} and the assertion on
order ideals of minimal generators of syzygies of modules of finite
projection dimension, introduced and proved by Evans and Griffith
(\cite{E-G1}, \cite{E-G2}) for equicharacteristic local rings. The
monomial conjecture asserts that given a local ring $R$ and a system
of parameters $x_1,\dots,x_n$ of $R$, $(x_1\dots x_n)^t\not\in
\(x_1^{t+1},\dots, x_n^{t+1}\)$. In \cite{Ho1} Hochster proved this
assertion for the equicharacteristic case and proposed it as a
conjecture for local rings of all characteristics. He also showed
(\cite{Ho1}) that this conjecture is equivalent to the direct summand
conjecture for module-finite extension of regular rings. From the
equational point of view both these conjectures claim that the
polynomial equation $(X_1\dots X_n)^t -
\sum\limits^n_{i=1}Y_iX_i^{t+1} = 0$ ($X_i,Y_js$ are variables) cannot
have a solution $\{x_1,\dots,x_n,\,y_1,\dots,y_n\}$ in any noetherian
ring, unless the height of $(x_1,\dots,x_n)$ is less than $n$. In the early
eighties in order to prove their syzygy
theorem (\cite{E-G1}), \cite{E-G2}) Evans and Griffith proved an
important aspect of order ideals of syzygies over equicharacteristic
local rings. We would like to generalize and state their result as a
conjecture on arbitrary local rings in the following way.

\proclaim{Order ideal conjecture}Let $(R,m)$ be a local ring. Let $M$
be a finitely generated module of finite projective dimension over $R$
and let $\Syz^i(M)$ denote its $i^{\text{th}}$ syzygy for $i>0$. If $\beta$
is a minimal generator of $\Syz^i(M)$, then the order ideal
$\ocal_{\Syz^i(M)}(\beta)(= \{f(\beta)|f\in \Hom_R(\Syz^i(M),R)\})$ has $\grade
\ge i$.
\endproclaim

For their proof of the syzygy theorem Evans and Griffith actually
needed a particular case of this conjecture: $M$ is locally free on the
punctured spectrum and $R$ is regular local. They reduced the proof of
the above conjecture for this special case to the validity of the
improved new intersection conjecture (also introduced by them) and
proved its validity over equicharacteristic local rings by using big
Cohen-Macaulay modules.
In the mid-eighties Hochster (\cite{Ho3}) proposed a new conjecture that is
deeply homological in nature---the canonical element conjecture. In
this conjecture Hochster assigns a canonical element $\eta_R$ to every
local ring $R$ and asserts that $\eta_R\ne 0$. He proved this
conjecture for the equicharacteristic case. In the same paper Hochster
(\cite{Ho2}) showed that the canonical element conjecture is equivalent to the
monomial conjecture and it implies the improved new intersection
conjecture due to Evans and Griffith. Later this author proved the
reverse implication (\cite{D1}).
%
%
%
Thus, the monomial conjecture implies a
special case of the order ideal conjecture: the case when $M$ is
locally free on the punctured spectrum of~$R$.

In this article we would prove that the validity of the order ideal
conjecture over regular local rings (actually a special case of it)
implies the monomial conjecture. First we would like to propose the
following definition.

Consider a finitely generated module $M$ of finite projective
dimension over a local ring $R$ with a minimal free resolution
$(F_\bullet,\,\beta_\bullet)$. Let $\Syz^i(M)$ denote the
$i^{\text{th}}$ syzygy of $M$ in $F_\bullet$. We say that for $i\ge 1$, $\Syz^i(M)$
satisfies the property ({\bf 0}) or simply ({\bf 0}) if for every
minimal generator $\alpha$ of $\Syz^i(M)$, the ideal $I$ generated by the
entries of $\alpha$ in $F_{i-1}$ has $\grade \ge i$. In our main theorem
(Theorem 1.4) we prove the following:

\proclaim{Theorem}The monomial conjecture is valid for all local rings
if for every almost complete intersection ideal $J$ of height $d$ in
any regular local ring $(R, m)$, $d<$ dimension of $R$, $\Syz^{d+1}(R/J)$
satisfies \rm({\bf 0}).
\endproclaim

The idea involved in the proof of this theorem evolved gradually through our work
in (\cite{D5}, \cite{D7}, \cite{D-G}). Actually our assertion is more specific than
what is stated above. We derive the following from the proof of the above theorem.

\proclaim{Corollary {\rm (Notations as above)}}The monomial conjecture is
valid for all local rings if for every almost complete intersection ideal
$J$ of height $d$ in a regular local ring $(R, m)$, $K_{d+1}(J; R)\bigotimes k \to \Tor_{d+1}^R(R/J, k)$
is the 0-map where $K_\bullet(J; R)$ is the Koszul complex corresponding to $J$ in $R$ and $k = R/m$.
\endproclaim

Next we prove a proposition (Prop. 1.5) that replaces the previous assertion
involving $J$ with the canonical module $\Omega$ of $R/J$. Let $J = (x_1, \dots, x_d, \lambda )$ and let
$\bx = (x_1, \dots, x_d)$ denote the ideal generated by an $R$-sequence. As a corollary we derive:

\it The monomial conjecture is valid if  $\Tor^R_d (k, R/\bx) \to \Tor^R_d (k, R/(\bx + \Omega))$ is the zero map
(equivalently $\Tor_d^R (k, \Omega) \to \Tor_d^R (k, R/\bx)$ is non zero).

\rm
The following proposition (prop. 1.6) provides a characterization of the above statement
in terms of syzygies of canonical modules.

\proclaim {Proposition {\rm (Notations as above)}} \it If $K_{d+1} (J; R) \bigotimes k \to \Tor_{d+1}^R (R/J, k)$
is the 0-map then $\Syz^d(\Omega)$ (minimal) has a free summand. Conversely, if $I$ is an ideal of $R$
of height d such that $R/I$ satisfies the Serre-condition $S_2$ and
its canonical module $\Omega$ is such that $\Syz^d(\Omega)$
has a free summand, then $\Tor_d^R (\Omega, k) \to \Tor_d^R (R/\bx ,\, k)$  is non-zero,
where $\bx$ denotes the ideal generated by a maximal $R$-sequence in $I$.
\endproclaim

As a corollary to our next theorem we observe the following:

\it Let (R, m) be a equicharacteristic regular local ring. Let
$I$ be an ideal in $R$ of codimension $d$ and let $\Omega$ denote the
canonical module of $R/I$. Then $\Syz^d(\Omega)$ has a free summand.

\proclaim {Remark} The above assertion is also valid in the graded equicharacteristic case
via the same mode of proof.
\endproclaim

In this connection let us recall that for any finitely generated module $M$ of
grade $g$, no $\Syz^i(M)$ can have a free summand for $i<g$ ([D2]).

Our next theorem describes the special cases where we can at present prove that $\Syz^d (\Omega)$
possesses a free summand.

\proclaim{Theorem} Let $(R, m, k)$ be a regular local ring in mixed characteristic $p > 0$ and let $I$ be an ideal
of height $d$ in $R$. Let $\bx$ denote the ideal generated by a maximal $R$-sequence contained
in $I$ and Let $\Omega = \Hom_R (R/I, R/\bx)$ denote the canonical module of $R/I$. Then $\Tor_d^R (\Omega, k) \to
\Tor_d^R (R/\bx, k)$ is non-zero (equivalently $\Syz^d (\Omega)$ possesses a free summand) in the following cases:
1) $\Omega$ is $S_3$ and 2) the mixed characteristic $p$ is a non-zero-divisor on $\Ext^{d+1} (R/I, R)$.
\endproclaim

In our last proposition (1.7) we derive a sufficient condition for
$\Syz^d(\Omega)$ to possess a free summand. We prove that \it if $\Omega^\prime$
denotes the lift of $\Omega$ in $R$ via $R \to R/\bx$ and if at least one of
$x_1,\dots,x_d$ is contained in $m\Omega^\prime$, then $\Syz^d(\Omega)$ has a free summand.

\rm
Throughout this work ``local'' means noetherian local. Over the years
different aspects of the monomial conjecture have been studied,
special cases have been proved and new equivalent forms have been
introduced (see \cite{Bh}, \cite{Br-H}, \cite{D1}, \cite{D2}, \cite{D3},
\cite{D4}, \cite{D5}, \cite{D6}, \cite{D-G}, \cite{Go}, \cite{He}, \cite{K},
\cite{O}, \cite{R3}, \cite{V}). Statements of four equivalent forms of
this set of conjectures, as proposed by Hochster in \cite{Ho3} are
given below for the convenience of the reader.


\subhead{A.\enspace Direct Summand Conjecture (DSC)}\endsubhead

Let $R$ be a regular local ring, and let $i:R\hra A$ be a module-finite
extension of $R$. Then $i$ splits as an $R$-module map.

\subhead{B.\enspace Canonical Element Conjecture (CEC)}\endsubhead

Let $A$ be a local ring of dimension $n$ with maximal ideal $m$ and
residue field $k$. Let $S_i$ denote the $i$th syzygy of $k$ in a
minimal resolution of $k$ over $A$, and let
$\theta_n:\Ext^n_A(k,\,S_n)\to H^n_m(S_n)$ denote the direct limit
map. Then $\theta_n$ (class of the identity map on $S_n$) $\ne 0$.

\subhead{C.\enspace Improved New Intersection Conjecture (INIC)}\endsubhead
Let $A$ be as before. Let $F_{\bullet}$ be a complex of finitely
generated free $A$-modules,
$$
F_{\bullet}:0\to F_s\to F_{s-1}\to \cdots\to F_1\to F_0\to 0,
$$
such that $\ell(H_i(F_{\bullet}))<\infty$ for $i>0$ and
$H_0(F_{\bullet})$ has a minimal generator annihilated by a power of
the maximal ideal $m$. Then $\dim A\le s$.

\subhead{D}\endsubhead
Let $(A,m,k=A/m)$ be a local ring of dimension $n$ which is a
homomorphic image of a Gorenstein ring. Let $\Omega$ denote the
canonical module of $A$. Then the direct limit map
$\Ext^n_A(k,\Omega)\to H^n_m(\Omega)$ is non-null.

Recall that for any equivalent version of this set of conjectures one
can assume that the local ring $A$ is a complete local (normal) domain. We have
$A = R/\tilde P$, where $R$ is a complete regular local ring. Let $S
= R/\underline{y}$, where $\underline{y}$ is the ideal generated by a
maximal $R$-sequence contained in $\tilde P$. Then $A=S/P$, $P =
\tilde P/\underline{y}$. Let $\Omega$ denote the canonical module for
$A$; $\Omega$ can be identified with $\Hom_S(A,S)$---an ideal of
$S$. Let $E$ denote the injective hull of the residue field of S (resp. R)
over S (resp. R). For any $S$-module ($R$-module) $T$, we write
$T^\vee$ to denote $\Hom_S(T,E)$ $(\Hom_R(T,E))$ and $\dim T$ to
denote the Krull dimension of $T$. These notations will be utilized
throughout this article.

\head{Section 1}\endhead

Our first proposition is an observation due to Strooker and
St\"{u}ckrad on a characterization of the monomial conjecture
(\cite{Str-St\"{u}}). This author independently proved a similar
characterization for the direct summand conjecture (\cite{D5}). Since
our main focus is on a proof of the monomial conjecture we provide a
proof of this result here.

\proclaim{1.1\enspace Proposition {\rm (Th. \cite{Str-St\"{u}})}}
With notations as above, $A$ satisfies MC if and only if for any
system of parameters $x_1,\dots,x_n$ of $S$, $\Omega$ is not
contained in $(x_1,\dots,x_n)$.
\endproclaim

\demo{Proof}Let $\xi_1,\dots,\xi_n$ be a system of parameters of $A$. We
can lift it to a system of parameters $x_1,\dots,x_n$ for $S$ such
that $\im(x_i) = \xi_i$, $1\le i\le n$. Conversely, any system of
parameters for $S$ is a system of parameters for $A$. Let us write
$\underline{x} = (x_1,\dots,x_n)$ and $\underline{\xi} =
(\xi_1,\dots,\xi_n)$. The monomial conjecture for $A$ is equivalent to the
assertion that, for every system of parameters $\xi_1,\dots,\xi_n$ of $A$,
the direct limit map $\alpha:A/\underline{\xi}\to H^n_m(A)$ is non-null
\cite{Ho1}. Because $S$ is a complete intersection, the direct limit
map $\beta:S/\underline{x}\to H^n_{m_S}(S)$ is non-null ($m_S$ = maximal ideal of $S$).
We have the following commutative diagram
$$
\matrix
S/\underline{x}\hfill\hfill & @>{\quad\beta\quad}>>\hfill\hfill &
H_{m_S}^n(S)\\
@VV{\overline \eta}V  @VV{\gamma}V\hfill\hfill\\
A/\underline{\xi}\hfill\hfill & @>{\quad\alpha\quad}>>\hfill\hfill &
H_m^n(A)\hfill\hfill\\
\noalign{\vskip-30pt}
& \searrow\\
\noalign{\vskip30pt}
\endmatrix
$$
where $\overline \eta$ is induced by the natural surjection $\eta : S \to A$
and $\gamma = H^n_{m_S}(\eta)$. This implies that $\alpha$ is non-null
$\Leftrightarrow \alpha\circ\eta = \gamma\circ\beta$ is
non-null $\Leftrightarrow H^n_m(A)^\vee\to (S/\underline{x})^\vee$
is non-null $\Leftrightarrow \im(\Omega\to S/\underline{x})$
is non-null $\Leftrightarrow \Omega\not\subset \underline{x}$
(recall that, by local duality, $H^n_m(A)^\vee = \Omega$).
\enddemo

In our next proposition we deduce the validity of the monomial conjecture
for all local rings from the validity of the same for all local almost
complete intersections. We prove the following:

\proclaim{1.2\enspace Proposition}
The monomial conjecture is valid for all local rings if and only if
it holds for all local almost complete intersections.
\endproclaim

\demo{Proof}Suppose MC holds for all local almost complete
intersections. Let $A$ be a complete local domain. Then we have $A =
R/\tilde P$, where $R$ is a complete regular local ring. We can choose
$y_1,\dots,y_r$---a maximal $R$ sequence in $P$ in such a way that
$\Tilde{P}R_{\tilde P} = (y_1,\dots,y_r)R_{\tilde P}$. Write $S =
R/\underline{y}$, where $\underline{y} = (y_1,\dots,y_r)$ and
$P = \tilde P/\underline{y}$. Then $S$ is a complete intersection, $A
= S/P$, $\dim S = \dim A$, and $PS_P = 0$. Let $\Omega = \Hom_S(A,S)$,
the canonical module of $A$. Consider the primary decomposition of (0)
in $S$: $0 = P\cap q_2\cap\cdots\cap q_h$, where $q_i$ is
$P_i$-primary and $\height \;P_i = \height\;P = 0$ for $2\le i\le h$. It can
be checked easily that $\Omega = q_2\cap\cdots\cap q_h$. Choose
$\lambda\in P - \bigcup_{i\ge 2}P_i$. Then $\Omega = \Hom(S/\lambda S,
S)$, and $S/\lambda S$ is an almost complete intersection. Since
$S/\lambda S$ satisfies MC by assumption, it follows from the above
proposition that $\Omega$ is not contained in the ideal generated by
any system of parameters in $S$. Hence, by Prop. 1.1, $A$ satisfies~MC.
\enddemo

Now we reduce the assertion in the monomial conjecture to a length
inequality between $\Tor_0$ and $\Tor_1$ on regular local rings.

\proclaim{1.3\enspace Proposition} The monomial conjecture is valid for all local rings if
and only if for every regular local ring $R$ and for every pair of
ideals $(I, J)$ of $R$ such that i)\thinspace $I$ is a complete
intersection, ii)\thinspace $J$ is an almost complete intersection
(i.e., $J$ is minimally generated by (height $J+1$) elements),
iii)\thinspace height $I$ + height $J = \dim R$ and iv)\thinspace
$(I+J)$ is primary to the maximal ideal of $R$, the following length
inequality holds:
$$
\ell\(R/(I+J)\) > \ell\(\Tor^R_1(R/I,\,R/J)\).
$$
\endproclaim

\demo{Proof}First assume that every pair of $(I, J)$, as in the
statement of our theorem, satisfies the length inequality:
$$
\ell(R/(I+J))>\ell(\Tor^R_1(R/I,R/J)).
$$
By Prop.\ 1.2, we can assume that $A$ is an almost complete intersection
ring of the form $S/\lambda S$, where $S$ is a complete intersection
and $\dim S = \dim A$. Let $\Omega = \Hom(S/\lambda S, S)$ denote the
canonical module for $A$. We consider the short exact sequence
$$
0\to S/\Omega @>f>> S\to S/\lambda S\to 0,
\tag{1}
$$
where $f(\overline{1}) = \lambda$. Let
$x_1^{\prime\prime},\dots,x_n^{\prime\prime}$ be a system of parameters
for $A$. We can lift $x_1^{\prime\prime},\dots,x_n^{\prime\prime}$ to
$x_1^\prime,\dots,x_n^\prime$ in $S$ in such a way that
$\{x_1^\prime,\dots,x_n^\prime\}$ form a system of parameters in
$S$. Let $\underline{x}^{\prime\prime} =
(x_1^{\prime\prime},\dots,x_n^{\prime\prime})$ and
$\underline{x}^\prime = (x_1^\prime,\dots,x_n^\prime)$. Tensoring (1)
with $S/\underline{x}^\prime$ yields the following exact sequence
$$
0\to \Tor_1^S(S/\underline{x}^\prime, S/\lambda S)\to
S/(\Omega+\underline{x}^\prime) @>{\overline f}>>
S/\underline{x}^\prime\to S/(\underline{x}^\prime + \lambda S)\to
0,\tag{2}
$$
where $\overline{f}$ is induced by $f$.
Then we have the following
$$
\Omega\not\subset \underline{x}^\prime \Leftrightarrow
\ell(S/(\underline{x}^\prime + \lambda S)) >
\ell(\Tor_1(S/\underline{x}^\prime, S/\lambda S)).
\tag{3}
$$

As in the proof of the previous proposition, let $R$ be a
complete regular local ring mapping onto $S$. Now lift
$x_1^\prime,\dots,x_n^\prime$ to an $R$-sequence
$x_1,\dots,x_n$ in $R$. Write $I = \underline{x}$ and $J =
(\underline{y}, \lambda)$. Then (3) translates to
$\ell(R/(I+J))>\ell(\Tor_1(R/I, R/J))$, as required in our statement.

For the converse part of the theorem, write $I = (x_1,\dots,x_n)$ and
$J = (y_1,\dots,y_r,\lambda)$, where $x_1,\dots,x_n$ and
$y_1,\dots,y_r$ form $R$-sequences such that $n+r = \dim R$. Let
$S=R/(y_1,\dots,y_r)$ and $A = S/\lambda S$, and let $x_i^\prime =
\im(x_i)$ in $S$ for $1\le i\le n$. Write $\underline{x}^\prime =
(x_1^\prime,\dots,x_n^\prime)$ and $\underline{x} =
(x_1,\dots,x_n)$. Since $I+J$ is primary to the maximal ideal and
since both $R/\underline{x}$ and $S$ are complete intersections, it
follows from \cite{Se} that $\Tor_i^R(R/\underline{x}, S) = 0$ for
$i>0$. This implies that $\Tor_i^R(R/I, R/J) =
\Tor_i^S(S/\underline{x}^\prime, A)$ for $i\ge 0$. Let $\Omega =
\Hom_S(A, S)$, the canonical module for $A$. Now (1)--(3) and the
subsequent arguments complete the proof.
\enddemo

\proclaim{1.4} \rm Our main theorem is the following.

\medskip

\endproclaim
%
%

\proclaim{Theorem}The monomial conjecture is valid for all local rings
if for every almost complete intersection ideal $J$ of height $d$ in
any regular local ring $R,\,d<\dim R,\,\Syz^{d+1} (R/J)$ satisfies~\rm({\bf 0}).
\endproclaim

\demo{Proof}By Prop.1.3 it is enough to prove the following.
\enddemo

\proclaim{Theorem}Let $(R,m)$ be a regular local ring and let $I,J$ be
two ideals of $R$ such that i)\th $R/I$ is Gorenstein and $J$ is an
almost complete intersection, ii)\th $I+J$ is $m$-primary and iii)\th
$\height I + \height J = \dim R$. Let $n = \height I$ and $d = \height J$.
Suppose that $d<\dim R$ and $\Syz^{d+1} (R/J)$ satisfies \rm({\bf 0}). Then
$$
\ell\(R/(I+J)\) > \ell\(\Tor_1^R\(R/I,\,R/J\)\).
$$
\endproclaim

\demo{Proof}Let $J = (y_1,\dots,y_d,\lambda)$ where $\{y_1,\dots,y_d\}$ form an $R$-sequence. Let
$(K_{\bullet},\gamma_{\bullet}) =
K_{\bullet}(y_1,\dots,y_d,\,\lambda;R)$ denote the Koszul complex
corresponding to $y_1,\dots, y_d,\,\lambda$ and let $H =
H_1(y_1,\dots,y_d,\lambda;R)$. Let $L_{\bullet}$ be a minimal free
resolution of $H$ and let $\psi_{\bullet}:L_{\bullet}\to K_{\bullet}(1)$ be a
lift of $H\hra \coker \gamma_2 = G$. We have the following
commutative diagram
$$
\vbox{\hbox{\vbox{\halign{\hbox to 28pt{}\eightpoint\tabskip=3pt\!
$#$\hfill\hfill & $#$\hfill\hfill & $#$\hfill\hfill & $#$\hfill\hfill &
$#$\hfill\hfill & $#$\hfill\hfill & $#$\hfill\hfill & $#$\hfill\hfill &
$#$\hfill\hfill & $#$\hfill\hfill & $#$\hfill\hfill & $#$\hfill\hfill &
$#$\hfill\hfill & $#$\hfill\hfill & $#$\hfill\hfill & $#$\hfill\hfill &
$#$\hfill\hfill & $#$\hfill\hfill\tabskip=0pt\cr
L_\bullet: & \longrightarrow & R^{r_{d+1}} &
\longrightarrow & R^{r_d} & \longrightarrow & R^{r_{d-1}} & \longrightarrow
 & \cdots & \longrightarrow & R^{r_0} & \longrightarrow & H
 &\longrightarrow & 0\cr
&&&& \Big\downarrow{\psi_d} && \Big\downarrow{\psi_{d-1}} &&&&
\Big\downarrow{\psi_0} && \Big\downarrow^{\hbox{\hskip-4.25pt}^{_\cap}}\cr
K_\bullet (1): && 0 & \longrightarrow & R & @>{\gamma_{d+1}}>> & R^{d+1}
& \longrightarrow & \cdots & @>{\gamma_2}>> & R^{d+1} &
\longrightarrow & G & \longrightarrow & 0\cr
\noalign{\vskip-68pt}
\noalign{\hfill\hfill\hfill
}\cr
\noalign{\vskip22pt}
}}}}
$$
\tenpoint\baselineskip17pt

\medskip

The mapping cone of $\psi_{\bullet}$ is a free resolution of $J$ and
thereby provides a free resolution of $R/J$ from which a minimal
resolution $(F_{\bullet},\beta_{\bullet})$ of $R/J$ can be extracted.
\enddemo

\proclaim{Claim}$\im \psi_d = R$.
\endproclaim

\demo{Proof of the Claim}If not, then $\im\psi_d\subset m$. Then
the copy of $R = K_{d+1}$ would be a free summand of $F_{d+1}$. Let
$\alpha$ be the image of $(1,\,\underline{0})$ in $\Syz^{d+1} (R/J) \subset
F_d$. Since $\beta_{d+1}(1,\,\underline{0})$ is a part or whole of
$\dsp \pmatrix y_1\\ y_2\\ \vdots\\ y_d\\ \lambda\endpmatrix$, the height of
the ideal generated by the entries of $\alpha$ must be less than
$(d+1)$. This contradicts our hypotheses that $\Syz^{d+1} (R/J)$ satisfies
({\bf 0}) and hence $\im\psi_d=R$.

This means that if $\eta_{\bullet}: K_{\bullet}\to F_{\bullet}$ lifts
the identity map on $R/J$, then $\eta_{d+1}(R)\subset m\,F_{d+1}$.

Recall that the spectral sequences
$\{\Ext^p_R(R/J,\,\Ext^q_R(R/I,R))\}$ and
$\{\Ext^i_R(\Tor^R_j(R/J,$\linebreak $R/I),\,R)\}$ converge to the
same limit when $p+q = i+j$. Since $R/I$ is Gorenstein, we have
$\Ext^i_R(R/I,\,R) = 0$ for $i\ne n$ and $\Ext^n_R(R/I,\,R)\simeq
R/I$. Since $\ell(R/(I+J))<\infty$ and $R$ is regular local of
dimension $(n+d)$, each of the above spectral sequence
degenerates. Hence
$$
\Ext_R^{n+d}\(\Tor_1^R(R/J,\,R/I),\,R\)\simeq \Ext_R^{d+1}(R/J,\,R/I).\tag{1}
$$
By (1) and local duality we have
$$
\ell\(\Tor_1^R(R/I,\,R/J)\)=\ell\(\Ext^{d+1}_R(R/J,\,R/I)\).
\tag{2}
$$

We also note that since $n+d=\dim R$ and $\ell\(R/(I+J)\)<\infty$,
$\Ext^i_R(H,\,R/I)=0$ for $i<d$.

The mapping cone of $\psi_{\bullet}$ leads to the following exact sequence:
$$
0\to K_{\bullet}\to
(\text{a~free~resolution~of~}R/J\text{~from~the~mapping~cone~of~}
\psi_{\bullet})\to L_{\bullet}(-2)\to 0.
$$
\noindent
Since $\Ext_R^{d-1}(H,\,R/I)=0$ applying $\Hom_R(-,\,R/I)$ to this
sequence we obtain the following short exact sequence
$$
0\to \Ext_R^{d+1}(R/J,\,R/I)\to H^{d+1}(y_1,\dots,y_d,\lambda;R/I) = R/(I+J).
$$
\noindent
This inclusion cannot be an isomorphism since $\eta_{d+1}(R)\subset m\,F_{d+1}$.

Thus, from (2), $\ell\(R/(I+J)\)>\ell\(\Tor_1^R(R/I,\,R/J)\)$ and the proof of
our theorem is complete.
\enddemo

We have the following corollary from the above proof.

\proclaim{Corollary}The monomial conjecture is valid if for every almost complete
intersection ideal $J$ of height $d$ in a regular local ring $(R, m, k = R/m)$,
$K_{d+1} (J; R) \bigotimes k \to \Tor_{d+1} (R/J, k)$ is the 0-map where $K_\bullet (J; R)$
is the Koszul complex corresponding to $J$ in $R$.
\endproclaim

The proof of this corollary follows from our observation after the diagram and
the claim in the proof of the above theorem.

\proclaim{1.5} \rm In our next proposition we study the assertion in the above corollary. \endproclaim
\proclaim{\enspace Proposition}Let $(R,m)$ be a local ring of $\dim
n$. Let $\{x_1,\dots,x_d\}\subset m$ be a regular sequence on $R$ and let
$\lambda$ be a $0$-divisor and not a parameter on $R/\bx$, $\bx=$
the ideal generated by $x_1,\dots,x_d$. Let
$\Omega=\Hom(R/(\bx,\lambda),\,R/\bx)$. Let
$K_\bullet(\bx,\lambda;R)=K_\bullet(x_1,\dots,x_d,\lambda;R)$
and $K_\bullet(\bx;R)$ denote the Koszul complexes corresponding
to $x_1,\dots,x_d,\lambda$ and $x_1,\dots,x_d$ respectively. Let
$(L_\bullet, c_\bullet)$ be a minimal free resolution of $\Omega$. Let
$\psi_\bullet:L_\bullet\to
K_\bullet(\bx,\lambda;\gamma_\bullet)(+1)$ and
$\phi_\bullet:L_\bullet\to K_\bullet(\bx;\delta_\bullet)$ be lifts
of $\Omega\simeq H_1(\bx,\lambda;R)\hra G=\Coker \gamma_2$ and
$\Omega\hra R/\bx$ respectively. Then $\im \psi_d=R$ if and only if
$\im \phi_d=R$.
\endproclaim

\demo{Proof}We have the following commutative diagrams:

\vglue-0.25truein

$$
\vbox{\hbox{\vbox{\halign{\hbox to 28pt{}\eightpoint\tabskip=3pt\!
$#$\hfill\hfill & $#$\hfill\hfill & $#$\hfill\hfill & $#$\hfill\hfill &
$#$\hfill\hfill & $#$\hfill\hfill & $#$\hfill\hfill & $#$\hfill\hfill &
$#$\hfill\hfill & $#$\hfill\hfill & $#$\hfill\hfill & $#$\hfill\hfill &
$#$\hfill\hfill & $#$\hfill\hfill & $#$\hfill\hfill & $#$\hfill\hfill &
$#$\hfill\hfill & $#$\hfill\hfill\tabskip=0pt\cr
L_\bullet: & \longrightarrow & R^{r_{d+1}} &
\longrightarrow & R^{r_d} & \longrightarrow & R^{r_{d-1}} & \longrightarrow
 & \cdots & \longrightarrow & R^{r_0} & \longrightarrow & H_1=\Omega
 & \longrightarrow & 0\cr
&&&& \Big\downarrow{\psi_d} && \Big\downarrow{\psi_{d-1}} &&&&
\Big\downarrow\phantom{{\Psi_0}} &&
\Big\downarrow^{\hbox{\hskip-4.25pt}^{_\cap}}\cr
\multispan{2}{$K_\bullet(\bx,\lambda;R)(+1):$} & 0 & \longrightarrow & R
& \longrightarrow & R^{d+1}
& \longrightarrow & \cdots & \longrightarrow & R^{d+1} &
\longrightarrow & G & \longrightarrow & 0\cr
\noalign{\vskip-36pt}
\noalign{\hfill\hfill\text{(1)}}\cr
}}}}
$$
\tenpoint\baselineskip17pt
\bigskip
\noindent and
$$
\vbox{\hbox{\vbox{\halign{\hbox to 28pt{}\eightpoint\tabskip=3pt\!
$#$\hfill\hfill & $#$\hfill\hfill & $#$\hfill\hfill & $#$\hfill\hfill &
$#$\hfill\hfill & $#$\hfill\hfill & $#$\hfill\hfill & $#$\hfill\hfill &
$#$\hfill\hfill & $#$\hfill\hfill & $#$\hfill\hfill & $#$\hfill\hfill &
$#$\hfill\hfill & $#$\hfill\hfill & $#$\hfill\hfill & $#$\hfill\hfill &
$#$\hfill\hfill & $#$\hfill\hfill\tabskip=0pt\cr
L_\bullet: & \longrightarrow & R^{r_{d+1}} &
\longrightarrow & R^{r_d} & \longrightarrow & R^{r_{d-1}} & \longrightarrow
 & \cdots & \longrightarrow & R^{r_0} & \longrightarrow & \Omega
 &\longrightarrow & 0\cr
&&&& \Big\downarrow{\phi_d} && \Big\downarrow{\phi_{d-1}} &&&&
\Big\downarrow\phantom{\Psi_0} &&
\Big\downarrow^{\hbox{\hskip-4.25pt}^{_\cap}}\cr
&& 0 & \longrightarrow & R & @>{\alpha_d}>> & R^d
& \longrightarrow & \cdots & \longrightarrow & R &
\longrightarrow & R/\bx & \longrightarrow & 0\cr
\noalign{\vskip-36pt}
\noalign{\hfill\hfill\text{(2)}}\cr
}}}}
$$
\tenpoint\baselineskip17pt

\bigskip

Suppose that $\im\psi_d=R$. Let $\eta_\bullet :
K_\bullet(\bx,\lambda;R)(+1)\to K_\bullet(\bx;R)$ denote the natural
surjection. Since $\eta_d=\Id_R$, $\eta_d\cdot \psi_d :
R^{r_d}\to R$ is surjective. Since $\eta_\bullet\cdot \psi_\bullet :
L_\bullet \to K_\bullet(\bx;R)$ lifts $\Omega\hra R/\bx$,
$\eta_\bullet\psi_\bullet$ is homotopic to $\phi_\bullet$. Hence
$\phi_d(R^{r_d})=R$. Conversely let $\im\phi_d=R$. Let
multiplication by $\lambda : K_\bullet(\bx;R)\to K_\bullet(\bx;R)$
lift the multiplication by $\lambda$ on $R/\bx$. Since
$\lambda\Omega=0$ in $R/\bx$, $\lambda_{\cdot}\phi_\bullet$ is
homotopic to the $0$-map: $L_\bullet\to K_\bullet(\bx;R)$. Let
$h_\bullet=(h_i)_{0\le i\le d-1}$ denote corresponding homotopies;
$h_i:R^{r_i}\to R^{\binom{d}{i+1}}$. Since
$K_\bullet(\bx,\lambda;R)$ is the mapping cone of
$\lambda:K_\bullet(\bx;R)\to K_\bullet(\bx;R)$, we can define
$\psi_\bullet^\prime : L_\bullet\to K_\bullet(\bx;\lambda;R)$, lifting
$\Omega=H_1\hra G$, in the following way:
$\psi_\bullet^\prime(R^{r_i})=(\phi_i,h_i)$. Then
$\psi_d^\prime(R^{r_d})=R$. Since
$\psi_\bullet^\prime,\psi_\bullet$ both lift the inclusion $\Omega\hra
G$, they are homotopic and hence $\psi_d(R^{r_d})=R$.
\enddemo

\proclaim{Corollary {\rm (with notations as above)}}Let $J = (x_1,\dots,x_d,\lambda)$
be an almost complete intersection ideal in a regular local ring $(R, m, k)$ and let $\Omega$
denote the canonical module of $R/J$. Then $K_{d+1}(J; R) \bigotimes k \to \Tor_{d+1}^R (R/J, k)$
is the 0-map if and only if $\Tor_d^R(R/\bx, k) \to \Tor_d^R(R/(\bx + \Omega), k)$ is the 0-map
or equivalently $\Tor_d^R (\Omega, k) \to \Tor_d^R (R/\bx, k)$ is non-zero.
\endproclaim

\proclaim {1.6} \rm The following proposition is partially a consequence of the above proposition.
\endproclaim

\proclaim {Proposition {\rm (Notations as in the previous corollary)}} \it If $K_{d+1} (J; R)
\bigotimes k \to \Tor_{d+1}^R (R/J, k)$ is the 0-map then $\Syz^d(\Omega)$ (minimal) has a free
summand. Conversely, if $I$ is an ideal of $R$ of height d such that $R/I$ is $S_2$ and
its canonical module $\Omega$ is such that
$\Syz^d(\Omega)$ has a free summand, then $\Tor_d^R (\Omega, k) \to \Tor_d^R (R/\bx ,\, k)$  is non-zero,
where $\bx$ denotes the ideal generated by a maximal $R$-sequence in $I$.
\endproclaim

\demo{Proof} Let us assume that $K_{d+1}(J; R) \bigotimes k \to \Tor^R_{d+1}(R/J, k)$ is the 0-map.
This implies, by diagram (1) in the proof of the above proposition, that
$\psi_d(L_d)=R$. Then we can write $L_d = R\bigoplus\/F$
where $F = \ker \psi_d$. Hence, from the commutativity of the diagram (1) in proposition (1.5),
it follows that $c_{d+1}(L_{d+1}) \subset\/F$ and this implies that $\Syz^d(\Omega)$
has a free summand.

Conversely let ($L_\bullet, c_\bullet)$ denote a minimal free resolution of $\Omega$.
If $\Syz^d(\Omega)$ has a free summand then $L_d$ has a free generator $e$ such that
$c_{d+1}^*(e^*) = 0$ ($-^* = \Hom_R(-,R)$). Let $L_\bullet^*$ denote the complex
$0 \to L_0^*\to L_1^*\to ... \to L_d^* \to G \to 0$. Since grade
$\Omega\geqslant\/d$, $L_\bullet^*$ is exact. Then im $e^*\in\Ext^d(\Omega, R)$
is non-zero and is a minimal generator of $G$. Let $x_1, \dots, x_d$ be a maximal $R$ sequence
contained in $I$ and let $\bx$ denote the ideal generated by them. Since $R/I$ is $S_2$, we have
$\Ext^d(\Omega, R) \cong \Hom (\Omega, \Omega) \cong R/I$. Let $G_{d+i} = \Coker c_{d+i}^*$ for
$i \geqslant 1$. We consider the following short exact sequences :
$$
0 \to R/I \to G \to \im c_{d+1}^* \to 0, \ \ 0 \to \im c_{d+1}^* \to L_{d+1}^* \to G_{d+1} \to 0,
$$
$$
0 \to \Ext^{d+1} (\Omega, R) \to G_{d+1} \to \im c_{d+2}^* \to 0, \ \ \dots,
$$
$$
0 \to \im c_{n-2}^* \to L_{n-2}^* \to \Ext^{n-2} (\Omega, R) \to 0.
$$
Since $\Omega$ is $S_2$, it follows from the above sequences that grade $\Ext^i (\Omega, R) \geqslant i+2$
for $i > d$ and hence it can be easily checked that
$\Omega = \Ext_R^d (G, R) \cong \Ext_R^d (R/I, R)$. Let $\alpha$ denote
the composite of $R/\bx \to R/I \cong \Ext_R^d(\Omega, R) \hookrightarrow G$, where im$1 \in R/\bx$ goes
to \linebreak im $e^* \in \Ext^d(\Omega, R)$. Let $\phi_\bullet^* : K_\bullet (\bx; R) \to L_\bullet^*$ denote
a lift of $\alpha$. Then $\phi_\bullet : L_\bullet \to K^\bullet (\bx; R)$ lifts the injection
$\Omega = \Ext^d (G, R) \hookrightarrow R/\bx$, where $\phi_d (L_d) = R$.
Thus, we obtain our required assertion.
\enddemo

In our next theorem we prove that for any ideal in a regular local ring the first part
of the above proposition is implied by the validity of the order ideal conjecture.

\proclaim{Theorem} Let us assume that the order ideal conjecture is valid for a regular
local ring $(R, m, k)$. Let $I$ be any ideal of $R$ of codimension d and let $\Omega$ denote the
canonical module of $R/I$. Then $\Syz^d(\Omega)$ has a free summand.
\endproclaim

\demo{Proof} Let $(F_\bullet, \phi_\bullet): \to F_d \to F_{d-1} \to \dots \to F_1 \to R \to 0$ denote a minimal
free resolution of $R/I$ over $R$. Let $F_\bullet^*$ denote the complex:
$0 \to R @>\; \phi_1^* \;>>F_1^* \to \dots \to F_d^* \to G_d \to 0$ where $G_{d+i} = \Coker \phi_{d+i}^*$,
for $0 \leq i \leq 1$. Since height of $I = d$, we have $H_i (F_\bullet^*) = 0$ for $i < d$ and $\Omega =
\Ext^d (R/I, R) \hookrightarrow G_d$. Let $L_\bullet$ be a minimal free resolution of $\Omega$ and let
$\psi_\bullet : L_\bullet \to F_\bullet^*$ lift $\Omega \hookrightarrow G_d$. Then $\psi_d (L_d) = R$;
for otherwise in the minimal free resolution $P_\bullet$ of $G_{d+1}$ obtained from the mapping cone of
$\psi_\bullet,  P_{d+1}$ would contain a copy of $R = L_0^*$ as a summand. Let $I = (y_1, \dots, y_s)$ where
$s = \ell (\Tor_1^R (k, R/I))$. Since height of $I = d$ and $\phi_1^* (1) = (y_1, \dots, y_s)$, this would
contradict the order ideal conjecture. Hence $\psi_d (L_d) = R$. This implies, by a previous
observation, that $\Syz^d (\Omega)$ has a free summand.
\enddemo

\proclaim {Corollary} Let $(R, m)$ be a equicharacteristic regular local ring and let I be an ideal
of R of codimension d. Let $\Omega$ denote the
canonical module of $R/I$. Then $\Syz^d(\Omega)$ has a free summand.
\endproclaim

\proclaim {Remark} The above assertion is also valid in the graded equicharacteristic case
via the same mode of proof.
\endproclaim

Our next theorem describe the cases where we are at present able to prove that $\Syz^d (\Omega)$ possesses a free summand.

\proclaim{Theorem} Let $(R, m, k)$ be a regular local ring in mixed characteristic $p > 0$ and let $I$ be an ideal
of height $d$ in $R$. Let $\bx$ denote the ideal generated by a maximal $R$-sequence contained
in $I$ and Let $\Omega = \Hom_R (R/I, R/\bx)$ denote the canonical module of $R/I$. Then $\Tor_d^R (\Omega, k) \to
\Tor_d^R (R/\bx, k)$ is non-zero (equivalently $\Syz^d (\Omega)$ possesses a free summand) in the following cases:
1) $\Omega$ is $S_3$ and 2) the mixed characteristic $p$ is a non-zero-divisor on $\Ext^{d+1} (R/I, R)$.
\endproclaim

\demo{Proof} Let $(F_\bullet, \phi_\bullet): \to F_d \to F_{d-1} \to \dots \to F_1 \to R \to 0$ denote a minimal
free resolution of $R/I$ over $R$. Let $F_\bullet^*$ denote the complex:
$0 \to R @>\; \phi_1^* \;>>F_1^* \to \dots \to F_d^* \to F_{d+1}^* \to F_{d+2}^* \to 0$.
Let $G_{d+i} = \Coker \phi_{d+i}^*$, for $0 \leq i \leq 2$. We have following short exact sequences for
$0 \leq i \leq 2$:
$$
0 \to \Ext^{d+i} (R/I, R) \to G_{d+i} \to \im\phi_{d+i+1}^* \to 0.
$$
Let $F_\bullet^\prime$ denote the complex $F_\bullet^*$ truncated at the $d$th spot i.e.
$H_0 (F_\bullet^\prime) = G_d$. Since height of $I = d, F_\bullet^\prime$ is a minimal free resolution of $G_d$.
Let $L_\bullet$ be a minimal free resolution of $\Omega$ and let $\psi_\bullet : L_\bullet \to F_\bullet^\prime$
denote a lift of the inclusion $\Omega = \Ext^d (R/I, R)\hookrightarrow G_d$.

{\bf {Claim.}} Suppose that either $\Omega$ is $S_3$ or $p$ is a non-zero-divisor on Ext$_R^{d+1} (R/I, R)$.
Then $\psi_d(L_d) = R$.
\enddemo

\demo{Proof of the Claim} If $\psi_d(L_d) \neq R$, then in the minimal free resolution $P_\bullet$ of $G_{d+1}$
extracted from the mapping cone of $\psi_\bullet$, $P_{d+1}$ would contain the copy of $R = F_0^\prime$ as a
free summand. Let $I = (y_1, \dots, y_{r_1})$ where $r_1$ = $\Tor_1^R (R/I, k)$. Since
height of $I = d$ and $\phi_1^* (1) = (y_1, \dots, y_{r_1})$, this would imply that $\Syz^{d+1} (G_{d+1})$
does not satisfy ({\bf 0}). Let $S = R/\bx$; then $\Omega = \Ext_R^d (R/I, R) = \Hom_S (R/I, S)$.
If $\Omega$ is $S_3$, then $\Ext_S^1 (R/I, S) = 0$, i.e. $\Ext^{d+1} (R/I, R) = 0$ and Hence $G_{d+1} \simeq$
Im $\phi_{d+2}^* \hookrightarrow F_{d+2}^*.$ Thus $p$ is a non-zero-divisor on $G_{d+1}$.
Let $\overline{R} = R/pR$ and $\overline{G_{d+1}} = G_{d+1}/pG_     {d+1}$. Since $\overline{G_{d+1}}$ has finite
projective dimension over $\overline{R}$ and $\overline{R}$ is equicharacteristic local ring,
$\Syz^{d+1} \overline{(G_{d+1})}$ must satisfy ({\bf 0}) by the order ideal theorem. This leads to a
contradiction and hence $\psi_d(L_d) = R$. If $p$ is a non-zero-divisor on $\Ext^{d+1}(R/I, R)$, then $p$
is a non-zero-divisor on $G_{d+1}$. Hence, arguing in a similar way as above, it follows that
$\psi_d (L_d) = R$ also in this case.

Let $\theta_\bullet$: $K_\bullet(\bx ; R) \to F_\bullet$ denote a lift of $R/\bx \to R/I$.
Dualizing $\theta_\bullet$ and combining the diagram corresponding to $\theta_\bullet^*$
with the same corresponding to $\psi_\bullet$ we get our required result.

\enddemo

\subhead{1.7}\endsubhead Finally we point out a sufficient condition for the assertion in the
above theorem to be valid in a more general set-up.

\proclaim{Proposition}Let $(R,m)$ be a local ring of dimension $n$ and
let $x_1,\dots, x_d\in m$ form an $R$-sequence. Let $\lambda$ be a
zero-divisor and not a parameter on $R/\bx$, $\bx=$ the ideal
generated by $x_1,\dots,x_d$ and let $\Omega=\Hom(R/(\bx,\lambda),\,R/\bx)$.
Let $\Omega^\prime$ be a lift of $\Omega$ in $R$ via the surjection $R\to R/\bx$. Let
$K_\bullet=K_\bullet(\bx;R;\beta_\bullet)$ and $(L_\bullet,\gamma_\bullet)$
denote the Koszul complex corresponding to $x_1,\dots,x_d$ and a minimal free resolution
of $\Omega$ respectively. Let $\phi_\bullet : L_\bullet \to K_\bullet$ denote a
lift of $\Omega\hra R/\bx$. If at least one of $x_1,\dots,x_d$ is
contained in $m\Omega^\prime$, then $\phi_d(L_d)=R$.
\endproclaim

\demo{Proof}Let $\overline{R}=R/x_1R$, $\overline{L}_\bullet =
L_\bullet\otimes \overline{R}$ and $\overline{K}_\bullet =
K_\bullet(\bx;R)\otimes \overline{R}$. Then
$H_1(\overline{L}_\bullet)=\Omega$, $H_i(\overline{L}_\bullet)=0$ for
$i>1$ and $H_1(\overline{K}_\bullet) = R/\bx = \overline{R}/\bx$
and $H_i(\overline{K}_\bullet)=0$ for $i>1$. We have a
canonical surjection $\eta_\bullet(=(\eta_i), i\ge 1) :\overline{K}_\bullet(+1)\to
K_\bullet(\overline{x}_2,\dots,\overline{x}_d;\overline{R})$ such that
$\eta_1$ induces the identity map on
$H_1(\overline{K}_\bullet)$ to $\overline{R}/\bx =
H_0(K_\bullet(\overline{x}_2,\dots,\overline{x}_d;\overline{R}))$. Let
$F_\bullet$ denote a minimal free resolution of $\Omega$ over
$\overline{R}$ and let $\theta_\bullet\(=(\theta_i)_{i\ge 0}\) :
F_\bullet\to \overline{L}_\bullet(+1)$ lift the natural inclusion
$\Omega = H_1(\overline{L}_\bullet)\hra \Coker
\overline{\gamma}_2$. Then $\eta_\bullet\overline{\phi}_\bullet
\theta_\bullet : F_\bullet\to
K_\bullet(\overline{x}_2,\dots,\overline{x}_d;R)$ lifts the natural
inclusion $\Omega\hra \overline{R}/\bx$ over $\overline{R}$. If
$\eta_d\overline{\phi}_d\theta_{d-1}(F_{d-1})=\overline{R}$, then it
follows that $\phi_d(L_d)=R$. Hence by induction on $d$ it is enough
to consider $d=1$ and $x=x_d$ is contained in $m\Omega^\prime$.

Let $\mu_1,\dots,\mu_h$ denote a minimal set of generators of
$\Omega^\prime$ over $R$. Since $x\in m\Omega^\prime$, there exists
$a_1,\dots,a_h\in m$ such that $x=\sum\limits^h_{i=1}a_i\mu_i$.

Consider the commutative diagram
$$
\vbox{\hbox{\vbox{\halign{\hbox to 18pt{}\eightpoint\tabskip=3pt\!
$#$\hfill\hfill & $#$\hfill\hfill & $#$\hfill\hfill & $#$\hfill\hfill &
$#$\hfill\hfill & $#$\hfill\hfill & $#$\hfill\hfill & $#$\hfill\hfill &
$#$\hfill\hfill & $#$\hfill\hfill & $#$\hfill\hfill & $#$\hfill\hfill &
$#$\hfill\hfill & $#$\hfill\hfill & $#$\hfill\hfill & $#$\hfill\hfill &
$#$\hfill\hfill & $#$\hfill\hfill\tabskip=0pt\cr
0 & \longrightarrow & \Syz^1 (\Omega) & \longrightarrow & R^h &
\longrightarrow & \Omega & \longrightarrow & 0\cr
& & \,\Big\downarrow{\phi_1} & & \,\Big\downarrow{\phi_0} & &
\Big\downarrow^{\hbox{\hskip-4.25pt}^{_\cap}}\cr
0 & \longrightarrow & R & @>\;x\;>> & R & \longrightarrow & R/xR &
\longrightarrow & 0\cr
}}}}
$$
where $\phi_0(e_i) = \mu_i$, $1\le i\le h$. Since $x=\sum{}a_i\mu_i$,
the element $\alpha=(a_1,\dots,a_h)\in \Syz^1 (\Omega)$ is such that
$\phi_0(\alpha)=x$. Hence $\phi_1(\Syz^1 (\Omega))=R$ and our proof is complete.
\enddemo



\enddocument